\documentstyle[12pt]{article}

\title{Quantum hyperboloid and braided modules}
\author{J.~Donin\\
Departement of Mathematics, Bar-Ilan University,\\
 52900 Ramat-Gan, Israel,\\
D.~Gurevich\\
ISTV, Universit\'e de Valenciennes\\
59304 Valenciennes, France,\\
V.~Rubtsov\\
Centre de Math\'ematiques, URA 169 du CNRS\\
 Ecole Polytechnique\\
91128 Palaiseau, France\\
and\\
ITEP, Bol.Tcheremushkinskaya 25,\\
 117259 Moscow, Russia}
\date{}

\begin{document}

\newtheorem{proposition}{Proposition}
\newtheorem{definition}{Definition}
\newtheorem{remark}{Remark}
\newcommand{\ren}{\rho_{End}}
\newcommand{\orr}{\overline{\rho}}
\newcommand{\gggg}{ g}
\newcommand{\tS}{\widetilde{S}}
\newcommand{\ahqc}{A_{h,q}^c}
\newcommand{\aqc}{A_{0,q}^c}
\newcommand{\ahq}{A_{h,q}}
\newcommand{\ahc}{A_{h,1}^c}
\newcommand{\ac}{A_{0,1}^c}
\newcommand{\aaa}{A_{0,1}^0}
\newcommand{\rn}{\rho_{\nu}}
\newcommand{\usl}{U_q(sl(2))}
\newcommand{\bC}{\mbox{$\bf C$}}
\newcommand{\bR}{\mbox{$\bf R$}}
\newcommand{\ug}{U(\gggg)}
\newcommand{\gug}{Gr\,U(\gggg)}
\newcommand{\us}{U(sl(2))}
\newcommand{\vv}{V^{\otimes 2}}
\newcommand{\uqs}{U_q(sl(2))}
\newcommand{\oI}{\overline{I}}
\newcommand{\uq}{U_q(\gggg)}
\newcommand{\sn}{sl(n)}
\newcommand{\ogg}{\overline{\gggg}}
\newcommand{\osl}{\overline{sl(2)}}
\newcommand{\aq}{A_{0,q}^0}
\def\ot{\otimes}
\def\De{\Delta}
\def\qq{(q+q^{-1})}

\maketitle

\begin{abstract}
When a quantum hyperboloid is realized, as a three - parameter
algebra $\ahqc$, in the usual manner, the following problem arises: what is a
"representation theory" of this algebra? We construct  the series of all spin
representations of $\ahqc$, and we discuss a braided version of the orbit
method, i.e. a  correspondence between orbits in $\gggg^*$ and $\gggg$-modules.
A braided trace and a braided involution are discussed as well.
\end{abstract}

\section{Introduction}

In the present paper we study a quantum hyperboloid from the point of view of
the generalized framework for  Quantum Mechanics suggested in \cite{GRZ}.
The main idea of the paper \cite{GRZ} was the following. When quantizing a
degenerate Poisson bracket we have to modify the ordinary
notions of Quantum Mechanics, namely, those of the Lie algebra, trace,
involution (conjugation) operator.

Such generalized objects and operators arise in
 particular as a result of the quantization
of some special Poisson pencils connected to the classical R-matrices, i.e.,
elements $R\in \wedge^2(\gggg)$ satisfying the classical Yang-Baxter equation
(CYBE)
$$[[R,R]]=[R^{12}, R^{13}]+[R^{12}, R^{23}]+[R^{13}, R^{23}]=0.$$
(There exist such generalized objects and operators which do not have any
quasiclassical nature, i.e., that can't be constructed by
quantzing some classical structures, cf. \cite{G}, but we will not discuss
them here.)

More precisely, given a representation $\rho :\gggg \to Vect(M)$
of a Lie algebra $\gggg$ in the space of vector fields on a manifold (variety)
$M$ equipped with a Poisson bracket $\{\,\,,\,\,\}$ which is preserved by the
fields of $Im(\rho)$, then the bracket
$$\{f,g\}_R=\mu<\rho^{\ot 2}(R), df\ot dg>, \,\, f,g\in Fun(M)$$
(where $\mu$ denotes the product in the algebra under consideration)
is also Poisson and compatible with the bracket $\{\,\,,\,\,\}$. If we want to
simultaneously quantize the Poisson pencil generated by these two brackets
we have to modify the scheme of the ordinary Quantum Mechanics.

Whereas in the symplectic case, there exists an invariant measure which
plays the role of a trace, while quantizing, in the case under consideration
such a measure does not exist.
This is the reason why we can't in general construct a trace with the usual
properties and have to use a twisted analogue of the
trace (and of other ingredients of Quantum Mechanics).

The principal aim of the present paper is to generalize this framework to the
case when the classical R-matrix entering the definition of the Poisson bracket
$\{\,\,,\,\,\}_R$ satisfies the modified CYBE, i.e., the elemenet
$[[R,R]]$ is $\gggg$-invariant. The simplest
quantum object arising from such Poisson bracket is a quantum hyperboloid.
More precisely it arises as a resuslt of the quantization on a usual
hyperboloid (regarded as an orbit in $\gggg^*$) of the Poisson pencil
generated by  the
Kirillov-Kostant-Souriau bracket $\{\,\,,\,\,\}_{KKS}$ and an R-matrix bracket
with $R=1/2(X\bigwedge Y)$ where the map $\rho$ is assumed to be
the coadjoint representation.

The final object of this procedure (often considered as a first stage of
quantization) is a three- parameter
algebra $\ahqc$ where $q$ is a braiding parameter, $h$ is a parameter of
quantization of the KKS bracket (it is introduced in the Lie bracket as a
factor) and $c$ labels the orbits (the case when $c=0$ corresponds to the
cone). More precisely, the algebra under consideration is a braided analogue
 of the quotient of the universal enveloping algebra $U(sl(2))$ by the ideal
generated by the element $C-c$, where $C$ is a Casimir element.

Let us note that such algebras were studied in a number of papers (\cite{P},
\cite{E}, \cite{NM}). These algebras, when equipped with some involution
operator, are called Podles' spheres (cf. Section 5 for a discussion of an
involution).

We are interested in the second stage of the quantization which usually
consists in an attempt to construct an
irreducible representation (or a series of those) of the quantum algebra.
Of course for any associative algebra, we can consider the left or right
regular representation but usually they are too big, i.e., reducible.

Concerning the algebras $\ahqc$, the following question is of great interest:
is it is possible to  construct their finite dimensional representations and
what are the possible values of the parameter $c$ for these representations?

In the classical case $(q=1)$ such a problem can be solved by means of
the representation theory of the  Lie algebra $sl(2)$.
Let us fix an $sl(2)$-module $U_k$ of spin $k$. Let $c_k$ be the corresponding
value of the Casimir element $C$. Then
the representation of the Lie algebra can be extended to a representation
of the quotient of $U(sl(2))$ by the ideal generated by
the element $C-c_k$.

In the present paper we construct a braided version of this approach.
First we describe a braided (or q-) analogue  of the Lie algebra $sl(2)$.
(It was introduced in \cite{DG1} and generalized to all Lie algebras in
\cite{DG2}). Then we construct a family of representations (of all spins)
of this Lie algebra-like object. The spaces $U_k$  where this braided Lie
algebra is represented are called {\em braided modules}.
These spaces are just $\usl$-modules but using the term "braided" we want to
emphasize the fact that they are equipped with the structure of a module over
a braided Lie algebra.

Let us emphasize that braided modules are suitable objects to develop a braided
version of the orbit method, which consists of a correspondence between
orbits in $\gggg^*$ and some $\gggg$-modules. We treat the algebra $\aqc$
with a fixed value of the parameter $c$ as a quantum or
braided orbit (namely, this algebra is a braided analogue  of a commutative
algebra,  cf. Section 2) and assign to such an "orbit" the braided module with
the same value of the Casimir operator. This is possible only for a series of
distinguished values, $c=c_k$, which are computed in Section 4.

Nevertheless we want to stress that the finitedimensional modules in the
classical case appear when one deals with
the $SU(2)$-orbits in $su(2)^*$. However, which a real form of the Lie algebra
and the Lie group is considered does not matter in the framework of the
algebraic approach, i.e., when one considers only algebraic functions on an
orbit and defines a trace as in the $su(2)$-case.

A braided deformation of the trace as well as of the ordinary involution are
considered in Section 5. If the construction of the braided trace is generally
accepted and not disputable it is not so for the  braided involution.
Usually one equips the algebra $\ahqc$ with an involution
 $*$ satisfying the classical property
$(ab)^*=b^*a^*$. Namely such a type of involution enters the construction
of Podles' spheres \cite{P}.

We state that for braided algebras $\ahqc$ it is more convenient to introduce
involutions with other properties which look like those of
super-algebras. Let us recall that in the latter case one usually imposes
the following axiom: $(ab)^*=(-1)^{\mid a \mid \mid
b \mid}b^*a^*$ where $\mid a \mid$ is the parity of $a$.

Moreover, if $V$ is a super-space then an involution $*$ acting in the space
$End(V)$ is compatible with the super-Lie bracket in the following sense
\begin{equation}
[\,\,,\,\,](*\ot *) =-*[\,\,,\,\,].
\end{equation}

We suggest a similar way to express the copatibility of an involution and
a braided structure in the algebra under consideration. Namely, we consider
the involutions compatible with the q-Lie braket in the sense of the
relation (1) and we classify them. As a corollary
we see that the compact form of the Lie algebra $sl(2, C)$ does not have
any braided deformation compatible with the braided structure in this sense.
(This , however, does not matter  for the computations of the braided trace.)

Remark that unlike the braided Lie algebra under
consideration, for the braided counterparts of other simple Lie algebras it
is not possible to construct a braided deformation of all finite dimensional
modules. It
seems very plausible that it is possible to do so only for such modules which
in the framework of the orbit method, correspond to symmetric orbits. (They
represent a particular case of so-called {\em R-matrix type orbits} in the
sense  of \cite{GP}.)

Let us remark also that the approach to the representation
theory of the algebra $\ahqc$ proposed in this paper enables us to consider
the construction of braided coadjoint vector fields suggested in \cite{DG1} and
\cite{DG2}  from a new point of view (see the remark at the end of Section 3.)

Throughout the paper the basic field $k$ is assumed to be
$R$ or $C$ unless specified otherwise.

\section{Quantum hyperboloid and braided commutativity}

  To construct a quantum hyperboloid it is sufficient to fix a representation
of the quantum group $\usl$ into a
threedimensional space $V$, decompose the space $\vv$ into a direct sum of
irreducible $\usl$-modules and impose a few natural equations on elements
of $\vv\oplus V\oplus k$ which are compatible with the action of the quantum
group $\usl$ and look like their classical counterparts.

Thus, let us consider the algebra $\usl$ generated by the three element $H, \,
X,\,Y$ satisfying the well-known relations
$$[H,X]=2X,\; [H,Y]=-2Y,\; [X,Y]=\frac{q^H-q^{-H}}{q-q^{-1}}$$
(through the paper the condition $q\not=0$ is assumed).
Let us equip this algebra with a copoduct defined on the basic elements
in the following way
$$\De(X)=X\ot 1+q^{-H}\ot X,\; \De(Y)=1\ot Y +Y\ot q^H,\;
\De(H)=H\ot 1+ 1\ot H.$$

It is well-known that this algebra has a Hopf structure, being equipped with
the antipode $\gamma$ defined by
$$\gamma(X)=-q^H X,\, \gamma(H)=-H,\, \gamma(Y)=-Yq^{-H}.$$

Let us consider a linear space $V$ with the base $\{u,\,v,\,w\}$ and  turn
$V$ into an $\usl$-module by setting
$$Hu=2u,\ Hv=0,\ Hw=-2w,\ Xu=0,\ Xv=-\qq u,\ Xw=v,$$
$$Yu=-v,\ Yv=\qq w,\ Yw=0$$
(it is easy to check that the above relations for $H, \,X,\,Y$ are
satisfied).

We want to stress that througthout this paper we deal with a
coordinate representation of module elements. We consider the endomorphisms
as matrices and their action as a left multiplication by these matrices.

Using the coproduct we can equip $\vv$ with a $\usl$-module structure as well.
This module is reducible and can be decomposed into a direct sum of three
irreducible $\usl$-modules
$$V_0={\rm span}((q^3+q)u w + v^2 + \qq w u),$$
$$V_1={\rm span}(q^2u v - v u,\ (q^3+q)( u w-w u) +(1-q^2)v^2,\
 -q^2v w + w v),$$
$$V_2={\rm span}(uu,\ uv+q^2vu,\  uw-q vv+q^4wu,\
vw+q^2 wv,\ ww)$$
of spins 0, 1  and 2 correspondinly (hereafter the sign $\ot$ is omitted).

Then only the following relations  imposed on the elements of the
space $\vv\oplus V\oplus k$ are coordinated with $\usl$-action:
$$C_q=(q^3+q)u w + v v + \qq w u=c,\,q^2u v - v u=-2 h u,$$
$$ (q^3+q)( u w-w u) +(1-q^2)v^2= 2hv,\,
 -q^2v w + w v=2hw$$
with arbitrary $h$ and $c$. The element $C_q$ will be called {\em braided
Casimir}.

Therefore it is natural to introduce a {\em quantum hyperboloid} as the
quotient algebra of the free tensor algebra $T(V)$ (if $h,\,q,\,c$ are fixed or
of the algebra  $T(V)[[h,q,c]]$ if $h,\,q,\,c$ are thought as formal
parameters)  over the ideal generated by elemens
$$(q^3+q)u w + v^2 + \qq wu-c,\,q^2u v - v u+2 h u,$$
$$ (q^3+q)( u w-w u)+(1-q^2)v^2-2hv,\, -q^2v w + w v-2hw.$$
This quotient algebra will be denoted by $\ahqc$.

Let us remark that the algebra $\ahc$ is a result of deformational quantization
of the KKS bracket, i.e., the restriction of the linear Poisson-Lie bracket
(given by the following multiplication table
$$\{v,\,u\}=2u,\, \{u,\,w\}=v,\,\{v,\,w\}=-2w)$$
on the variety $M$ fixed by $4uw+v^2=c$.
This variety is a manifold (namely, hyperboloid) when $c\not=0$ and
the cone if $c=0$.

The function algebra on these varieties (defined as restrictions of the
polynomials on $V$ to $M$) is just $\ac$. This algebra is commutative. A
natural
question arises: what is a braided (or q-)analogue  of this property?
(It is well-known that in symmetric tensor categories, i.e., those equipped
with
an involutive "transposition" $S$ a natural analogue  of it is a
{\em S-commutativity} condition expressed by the formula $\mu S=\mu$.)

In the case under consideration this property can be formulated by means of an
involutive operator $\tS:(\ahqc)^{\ot 2}\to (\ahqc)^{\ot 2}$. In fact this
operator is well-defined for any two $\usl$-modules $U$ and $V$ and sends $U\ot
V$ to $V\ot U$. It commutes with  action of $\usl$. We refer the reader to
\cite{DG2} for details, where the construction of similar operators is
discussed for any simple Lie algebra.

Since the algebra $\ahqc$ (for generic $q$) can be decomposed into a direct sum
of irreducible  $\usl$-modules  (what follows from flatness of
the deformation discussed below) the operator $\tS$ can be defined on the whole
$({\ahqc})^{\ot 2}$.

Now we are able to formulate an important property of the algebra $\aqc$: it
is $\tS$-commutative, i.e., $\mu \tS=\mu$ (we call this property {\em braided
commutativity}). This fact can be deduced from \cite{DS} (its rigorous
demonstration will be given elsewhere in more general context). Thus, roughly
speaking in the family $\ahqc$ only the algebra $\aqc$ can be treated as a
quantum hyberboloid (or cone) if by this we mean a  q-analogue  of the
commutative algebra $\ac$ (nevertheless according to tradition we use this
term in general case). Hence, we can treat the algebra $\aqc$ as a "classical"
object in the category of $\usl$-modules and $\ahqc$ as its quantum
counterpart.

Completing this section we want to make some remarks on the proof  of
flatness of deformation $\ac\to\ahqc$.
A proof given in \cite{DG1} was based on the following statements.

1. The algebra $\aq$ is Koszul (cf. \cite{BG} for definition). This fact was
proved in \cite{DG1} "by hand". Now there exists (for the case $q=1$ and
hence for a generic $q$) a more efficient proof valid for any simple Lie
algebra (cf. \cite{Be}, \cite{Bo}).

2. It is possible to describe the algebra $\ahqc$ as the enveloping algebra
of a generalized Lie algebra in the following sense (a little bit different
from that of \cite{DG2}). Let us consider the space $I= V_1 \oplus V_0$ and
introduce two maps $\alpha: I\to V$ and $\beta: I\to k$ as follows $\alpha:
V_0\to 0,\,\beta: V_1\to 0,$
 $$\alpha(q^2u v - v u)=-2 h u,\,\alpha((q^3+q)( u w-w u) +(1-q^2)v^2)=
2hv,\,$$ $$\alpha(-q^2v w + w v)=2hw,\, \beta ((q^3+q)u w + v^2 + \qq w u)=c.$$

We say that the data $(V,\, I\subset \vv,\, \alpha,\, \beta)$ define a
generalized Lie structure if the following relations are satisfied\\
\rm a. $Im(\alpha\ot
id- id \ot \alpha)(I \ot V\bigcap V\ot I) \subset I,$\\
\rm b. $(\alpha(\alpha\ot id- id \ot \alpha)+
\beta \ot id - id \ot \beta)(I \ot V\bigcap V\ot I)=0,$\\
\rm c. $\beta (\alpha\ot id- id \ot \alpha)(I \ot V\bigcap V\ot I)=0.$

Then in virtue of the main resul of \cite{BG} (which is in fact a slight
generalization of the PBW theorem in the form of \cite{PP}) we can deduce that
its graded adjoint algebra  $Gr\,\ahqc$ is isomorphic to $\aqc$.

Let us remark that the above conditions \rm a.-\rm c. represent the most
general analogue of the Jacobi  identity related to deformation theory.
However, they are useless from the representation theory point of view. On
 contrast, the Jacobi
identity  presented in the next Section is related to the representations of
the braided Lie algebra under consideration.

\section{q-Lie bracket and braided modules}

To prove the  flatness of the deformation discussed above we do not need that
 maps $\alpha, \beta$ be defined on all of $\vv$ but only on $I$.
Nevertheless there  exists a natural way to extend them to $\vv$.
We are only interested in the extension of the map $\alpha$ which will be
denoted
by $[\,\,,\,\,]$ and will be called {\em braided (or q-) Lie bracket}.

Setting $[\,\,,\,\,]=\alpha$ on $I$ and $[\,\,,\,\,]=0$ on $V_2$ we get a
morphism in the category of $\usl$-modules.

Let us remark that a similar q-Lie bracket can be defined for any simple Lie
algebra $\gggg$ (cf. \cite{DG2}).  It is possible as well to define the
corresponding  envelopping algebras in a natural way. Nevertheless it is a flat
deformation of its classical counterpart only for $\gggg=sl(2)$. If it is the
case the enveloping algebra denoted by $\ahq$ is defined in a similar way as
the
algebra $\ahqc$ but without the element $C_q-c$ between the generators of the
ideal.

Let us reproduce from \cite{DG1} the multiplication table for this bracket in
the base $\{u,\,v,\,w\}$:
$$
[u,u]=0,\ [u,v]=-q^2Mu,\ [u,w]=\qq^{-1}Mv,$$
$$[v,u]=Mu,\ [v,v]=(1-q^2)Mv,\ [v,w]=-q^2Mw,$$
$$[w,u]=-\qq^{-1}Mv,\ [w,v]=Mw,\ [w,w]=0,$$
where $M=2h(1+q^4)^{-1}$.

The space $V$ equipped with this bracket will be called {\em braided (or q-)
Lie
algebra} and denoted by $\osl$.

Let us note that this q-Lie algebra turns into the Lie algebra $sl(2)$ for
$q=1,\,\,h=2$ (and hence $M=2$).

Let us assign to any element $z\in V$ the "left adjoint" operator:
$\rho(z)x=[z,x]$. Then a natural question arises: whether the map $z\in V \to
\rho(z)\in End(V)$  defines a representation of the q-Lie algebra $\osl$, i.e.,
the relations $$q^2\rho(u)\rho(v)-\rho(v)\rho(u)=-2h\rho(u),\,
(q^3+q)(\rho(u)\rho(w)-\rho(w)\rho(u))+$$
$$
(1-q^2)\rho(v)^2=2h\rho(v),\,\,
-q^2\rho(v)\rho(w)+\rho(w)\rho(v)=2h\rho(w)$$
are satisfied?

The answer is negative. However,
it is possible to check  by straightforward computations
that this map defines for a generic $q$ a left {\em almost representation} of
this q-Lie algebra $\osl$ in the following sense.

\begin{definition}
We say that a map $\rho: V \to End(U)$ where $U$ is a $\usl$-module is
a almost representation of  the q-Lie algebra $\osl$ if it is
$\usl$-morphism  and there exists a factor $\nu\not=0$ such that
$$q^2\rho(u)\rho(v) -
\rho(v)\rho(u)=\nu(-2 h\rho(u)),\, (q^3+q)(\rho(u)\rho(w)-\rho(w)\rho(u))+$$
$$ (1-q^2)\rho(v)^2=\nu2h\rho(v),
-q^2\rho(v)\rho(w) + \rho(w)\rho(v)=\nu2h\rho(w)$$
(the case $\nu=1$ corresponds to a representation).
\end{definition}

\begin{remark} The fact that the above bracket defines an almost representation
means that the relations
$$ q^2[u,[v,z]]-[v,[u,z]]=-2\nu h[u,z],\, (q^3+q)([u,[w,z]]-[w,[u,z]]+
$$
$$(1-q^2)[v,[v,z]]=2\nu h [v,z],\,
-q^2[v,[w,z]]+[w,[v,z]]=2\nu h[w,z]$$
are satisfied for any $z$ and some $\nu$. This is another analogue of Jacobi
identity valid for the braided Lie algebra $\osl$. Note that for
Lie algebras $sl(n),\,\,n>2$ the orbit corresponding to the adjoint
representation in the frame of orbits method is not symmetric.
This is reason why it is hopeless to try to get a similar braided version
of Jacobi identity for these Lie algebras (cf. remark at the end of the
Introduction). \end{remark}

It is evident that the map ${\nu}^{-1}\rho$ is a representation of the
q-Lie algebra under question. Thus, extending the latter representation of
q-Lie algebra up to the whole algebra $\ahq$   we can construct a spin 1
representation of the algebra $\ahq$. By means of the above bracket it is also
possible to construct for a generic $q$ all integer spin almost
representations (and  consequently representations) of the algebra $\ahq$.

However, we will use here another, more direct, method enabling us
to construct the set of representations of all spins.

Let us fix a left spin $k$ irreducible  $\usl$-module $U=U_k$
and consider  the space $End(U)$ of endomorphisms of $U$ as
an $\usl$-module. This means that if $\rho : \usl \to End(U)$ is a
representation of the quantum group $\usl$ then $\ren :\usl \to End(End(U))$ is
defined as follows $$\ren(a)M=\rho(a_1)\circ
M\circ\rho(\gamma(a_2)),\,a\in\usl,
M\in End(U)$$ where $\circ$ denotes the matrix product, $\gamma$ is the
antipode in $\usl$ and $a_1\ot a_2$ is the Sweedler's notation for $\De(a)$.

Let us remark that this way to equip $End(U)$ with a structure of $\usl$-module
is coordinated with the matrix product in it, i.e.,
$$\ren(a)(M_1\circ M_2)=\ren(a_1)M_1\circ\ren(a_2)M_2.$$

Let us give the explicit form of the representation $\ren$:
$$\ren(X) M=\rho(X)\circ M-\rho(q^{-H})\circ M\circ\rho(q^{H})\circ\rho(X),$$
$$\ren(H) M=\rho(H)\circ M-M\circ\rho(H),\,
\ren(Y)M=(\rho(Y)\circ M-M\circ\rho(Y))\circ\rho(q^{-H}).$$

Let us decompose the $\usl$-module $End(U)$ into a direct sum of irreducible
$\usl$-modules. It is evident that for any spin $k$ in this decomposition there
is only an unique module isomorphic to $V$. We will call this condition {\em
unicity one}.

Let us define an $\usl$-morphism $\orr : V \to End(U)$ in a natural way
sending $V$ in the mentioned component of $End(U)$ (this
morphism is defined up to a factor).
\begin{proposition} The map $\orr$ is almost representation (for a generic
$q$).
\end{proposition}

{\bf Proof.}  By construction $\orr$ is an $\usl$-morphism. It is evident that
the elements
$$q^2\orr(u)\orr(v)-\orr(v)\orr(u),\,(q^3+q)(\orr(u)\orr(w)-\orr(w)\orr(u))+
(1-q^2)\orr(v)^2,$$
$$ -q^2\orr(v)\orr(w)+\orr(w)\orr(v) \in End(U)$$
generate an $\usl$-module  isomorphic to $V$ and therefore they coincide
correspondingly with $-\orr(u),\, \orr(v),\,\orr(w)$ up to a factor $\theta$ in
virtue of the unicity condition (we put $\theta=2h\nu$). The property that
$\theta\not=0$ for generic $q$ follows from the fact that it is valid for
$q=1$.
This complete the proof.

By the above method we can construct any spin  representation of the q-Lie
algebra under question and by extension we obtain all finite dimensional
representations of the  algebra $\ahq$.

\begin{remark}
Let us say a few words on the construction of "braided coadjoint vector fields"
suggested in \cite{DG1} and \cite{DG2}. It is well known that in the classical
case  a vector field is defined by means of Leibniz rule which can be
formulated
in terms of a coproduct operator defined on vector fields.
In the braided case, i.e., in the algebra $\ahq$ (unlike that $\usl$),
such coproduct does not
exist. So, even for a quantum hyperboloid if we define braided coadjoint vector
fields on  linear functions by means of the above q-Lie bracket it is not clear
what is a natural way to extend them on the higher power elements.

A method to do it "by analogy with the classical case" was suggested
in \cite{DG1},  \cite{DG2}.
In fact we can treat any braided coadjoint vector field as a
direct sum of almost representations of the braided Lie algebra in
finite dimensional modules. Now using the approach proposed above we can
coordinate the factors $\nu$ participating in these almost representations
(making them equal to the first one). This
method to coordinate the action of a braided coadjoint vector field on all
components of the "function space" on a quantum hyperboloid replaces
the mentioned above coproduct operator.
\end{remark}

\section{Braided Casimir}

In this Section we will generalize the
well-known property of the Casimir element stating that its image is a scalar
operator in any irreducible $sl(2)$-module to the braided case  and
compute the corresponding values of the braided Casimir.

Let $\orr=\orr_k$ denotes now the representation of the algebra $\ahq$
in the spin k module $U_k$.

\begin{proposition} Let $\orr : \ahq\to End(U_k)$ be such a representation.
Then the image $\orr(C_q)$ of the braided Casimir is a scalar operator for
generic $q$ (more precisely, it is so for any $q$ for which it is defined).
\end{proposition}

{\bf Proof.} Since $\orr$ is $\usl$-morphism and $C_q$ generates the trivial
$\usl$-module we have $$\orr(\rho(a)C_q))=\ren(a)\orr(C_q)=0\,\,\, {\rm
for\,\,\, any}\,\,\, a\in\usl$$
where $\rho : \usl\to End(U_k)$ is a representation of the quantum group
$\usl$.   For generic $q$ the elements $\ren(a) ,\,a\in\usl$ generate
the whole algebra $End(U_k)$. Using the above
explicit form of the representation $\ren$ it is easy to see that
$\orr(C_q)$ commutes with all elements of $End(U_k)$. This yields the
statement.

It is well-known that $dim\, U_k=l+1$ where $l=2k$.
Let us introduce some notations.
Denote by  $diag(a_1,\,a_2,...,a_{l+1})$
the diagonal matrices and by $diag_{\epsilon}(a_1,\,a_2,...,a_l),\,\,
\epsilon=\pm$ over-diagonal (if $\epsilon=+$) and sub-diagonal (if
$\epsilon=-$)
matrices.

Let us fixe the base in the $\usl$-module $U_k$
such that the corresponding representation $\rho=\rho_k:\usl\to End(U_k)$ is
of the form
$$\rho(X)=diag_+(1,\,1,\,...,1),\, \rho(H)=diag(l,\,l-2,\,...,-l),$$
$$\rho(Y)=diag_-(y_1,\,y_2,...,y_l)$$
where $y_1$ can be found from the following system
$$y_1=b_l,\, y_2-y_1=b_{l-2},..., y_l-y_{l-1}=b_{-l+2},\,-y_l=b_{-l},\,\,
b_i=(q^i-q^{-i})(q-q^{-1})^{-1}.$$

It is easy to check that the matrix $U=diag_+(q^{2(l-1)},\,q^{2(l-2)},...,1)$
satisfies the following conditions $\ren(X)U=0$ and $\ren(H)U=2 U$.
Let us consider the matrices  $V$ and $W$ such that
$-V=\ren(Y)U,\, \qq W=\ren(Y)U$. One can see that
$V=diag(v_1,\,v_2,...,v_{l+1})$
and $W=diag_-(w_1,...,w_l)$.

Using the  explicit form of the representation $\ren$ given above it is
possible
to find the values of all $v_i$ and $w_i$ but we need only  those $v_1$ and
$v_2$. We have  $v_1=y_1q^{l-2},\,\,v_2=y_2q^{l-2}-y_1q^l$.

Since the map $u\to U,\, v\to V,\, w \to W$ defines an almost
representation by virtue of the unicity
property the relations $$q^2UV-VU=-\theta U,\,(q^3+q)(UW-WU)+(1-q^2)V^2=\theta
V,$$ $$-q^2WV-WV=\theta W$$
are satisfied with some $\theta$. Substituting $U$ and $V$ to the first
relation and computing the first non-trivial matrix element we have
$$\theta=v_1-q^2v_2=y_1q^{l-2}-q^2(y_2q^{l-2}-y_1q^l)=y_1(q^{l+2}+q^{l-2})
-y_2q^l=q^{2l+1}+q^{-1}.$$
By the same reason we get
$(q^3+q)u_1 w_1+(1-q^2)v_1^2=\theta v_1$  using the second relation.

Therefore the first matrix element of the scalar operator
$(q^3+q)UW+V^2+\qq WU$ is equal to
$$(q^3+q)u_1 w_1+v_1^2=\theta v_1+q^2v_1^2=y_1q^{l-2}(\theta+y_1q^l)
=b_lb_{l+2}q^{2l-2}.$$

Thus, we have that the image of the braided Casimir under the above
almost representation is equal to $b_lb_{l+2}q^{2l-2}\,Id$.
We obtaine a representation of the braided Lie algebra $\osl$ if we put
$$\ren(u)=2h{\theta}^{-1}U,\,\ren(v)=2h{\theta}^{-1}V,\,
\ren(w)=2h{\theta}^{-1}W.$$

This provides the following
\begin{proposition} The value of the braided Casimir $C_q$ corresponding to the
braided $\ahq$-module $U_k$ is equal to
$$c_k=b_lb_{l+2}q^{2l-2}(2h{\theta}^{-1})^2,\,\,{\rm where}\,\,
\theta=q^{2l+1}+q^{-1},\, l=2k.$$
\end{proposition}

Let us consider two examples. If $l=1$ the matrices $U,\,V$ and $W$ are
correspondingly equal to
$$
\left(\begin{array}{cc}0&1\\0&0\end{array}\right),\,
\left(\begin{array}{cc}q^{-1}&0\\0&-q\end{array}\right),\,
\left(\begin{array}{cc}0&0\\q^{-1}&0\end{array}\right).$$
If $l=2$, they are
 $$\left(\begin{array}{ccc}0&q^2&0\\0&0&1\\0&0&0\end{array}\right),\,
\qq\left(\begin{array}{ccc}1&0&0\\0&1-q^2&0\\0&0&-q^2\end{array}\right),\,
\left(\begin{array}{ccc}0&0&0\\ 1&0&0\\0&1&0\end{array}\right).$$

\section{Discussion: braided trace and involution}

Let us remark that non-braided algebra $\ahc$ is multiplicity free, i.e.,
multiplicity of any $sl(2)$-module in it is not more than 1 (in fact only the
integer spin modules "live" in this algebra). A similar property  is true for
the algebra $\ahqc$. So there exists the unique (up to a factor) way compatible
with  $\usl$-action to introduce a {\em braided trace}  in this algebra as a
non-trivial operator $\ahqc\to k$ killing all  $\usl$-modules apart from the
trivial one.

This operator (denoted by $tr_q$) is a braided analogue  of the integral over
the
symplectic measure on a sphere  if we consider the "classical" braided algebra
$\aqc$ (recall that this algebra is $\tS$-comutative). It is a braided analogue
of the "quantum" trace if we consider the algebra $\ahqc$ or its
representations. We denote this operator by $Tr_q$ in the algebra $\ahqc$ and
$Tr_q^k$ in the $\ahqc$-module $U_k$.

Using a method of \cite{NM} it is easy to get the following

\begin{proposition} In the algebra $\aqc$ one has
$$tr_q v^{m}={{q^2-1}\over{2(q^{2m+2}-1})} (1+(-1)^m)(q{\sqrt c})^{-m}\, tr_q
1.$$ \end{proposition}

Applying a small modification of the method from \cite{NM} it is also possible
to
get a formula for  $Tr_qv^m$ but it is very more complicated and we do not
reproduce it. It would be interesting to get $Tr_q^k$ by direct computations
using the well-known definition of traces in braided categories and get a
braided version of the character formula.

\begin{remark} It is curious to note that classical character formula
can be obtained (at least for $sl(2)$-case) in a suitable form
as follows. We get in the classical case $(q=1)$ from the latter proposition
$tr e^{tv}=(e^{t\sqrt c}-e^{-t\sqrt c})(2t{\sqrt c})^{-1}\, tr 1$.

As for the trace $Tr^k e^{tV}$ where $V=diag(l,\,...,-l),\,\,l=2k$
(we put $h=2$) we have $Tr^k e^{tV}=(e^{(l+1)t}-e^{-(l+1)t})
(e^{t}-e^{-t})^{-1}$. Then the fraction $Tr^ke^{tV}/tr^ke^{tv}$ does not depend
on spin if we set $c=(l+1)^2=(2k+1)^2$ and assume that
$tr 1$ is proportional to $\sqrt c$. If it is the case we get
a formula $Tr^ke^{tV}=f(t)\, tr e^{tv}$ (with a factor $f(t)$ independed on
spin)
which looks like the character formula for $su(2)$. Taking into account that
the
value $\overline c$ of the Casimir $2(UW+WU)+V^2$ corresponding to the
$sl(2)$-module $U_k$ is $4k(k+1)=c-1$ we get a shift between values of Casimir.
In other words, the hyperboloid fixed by the relation $C=c$ "corresponds"
to the spin k $sl(2)$-module if the condition $\overline c=c-1$ is satisfied.

This correspondence is motivated by a version of the character formula and is
rather a "collective" phenomenon (our method does not allow us to get the above
correspondence if we consider an single orbit).
\end{remark}

We complete this Section with a discussion on involution operators in the
algebras under consideration. As it was mentioned in the Introduction one
usually  defines this operator by means of the "classical" property. Sometimes
one introduces an involution on a "quantum homogeneous space" by means of the
reduction of an involution operator defined in the quantum group or their
(restricted) dual object. As for the involutions defined in the latter objects
they are usually assumed to be coordinated with the Hopf structure.

We will consider here the involution operators which are coordinated with
the braided  structure in another way.
However, first we will recall the construction taken from \cite{GRZ} of
involutions  in the algebras living in a  tensor symmetric category, i.e., that
equipped with a involutive  "commutativity morphism" $S$.

Let $V$ be an object of such a category. Assuming the category to be rigid
let us consider the space $End(V)$ and involution
$*: End(V)\to End(V)$ in it satisfying two axioms
$*\mu=\mu (*\ot *)S$ and $S(id\ot *)= (*\ot id)S$. Then taking into account
that the "S-Lie bracket" $[\,\,,\,\,]$ is defined in the space $End(V)$
by $\mu(id-S)$ we can get the relation (1). (More precisely, we assume that
the space $V$ and all operators are given over the field $R$, then we consider
their complexification, cf. below.)

Let us remark that the relation (1) is universal: it has the same form in
the classical, the super- cases and in any symmetric category. Possessing now
a braided analogue  of Lie bracket we can extend this relation to the braided
case.

Let us consider
the space $\gggg=V$ over the field $k=R$ and assume the parameters $h$ and $q$
to be real. This means that all structure constants in the multiplication table
for the q-Lie bracket and all matrix elements of the morphisms
$S$ and $\tS$ (although we do not use them explicitly) are real as well.
Let us extend this bracket to the space $V_C=V\ot C$ by linearity. Let $*:
V_C\to
V_C$ be an involution (conjugation), i.e., an involutive $(*^2=id)$ operator
such that $(\lambda z)^*=\overline{\lambda}z^* ,\,\lambda\in C,\, z\in V_C$.
\begin{definition} We say that an involution $*$ is coordinated with the q-Lie
bracket $[\,\,,\,\,]$  if the relation (1) is satisfied.
\end{definition}

\begin{proposition} The odd elements with respect to this involution
(i.e., such that $z^*=-z$) form a subalgebra, i.e., the element $[a,\,b]$
is odd if $a$ and $b$ are.
\end{proposition}
{\bf Proof} is evident.

\begin{remark} One often considers involutions which differ from ours by the
sign. For such a type involution we have to change the sign in the relation (1)
and consider even elements instead of odd ones in the latter Proposition.
\end{remark}

Now we will classify all involutions $*:V_C\to V_C$ coordinated with q-Lie
bracket.

\begin{proposition} For a real $q\not=1$ there exist only two involutions
 in the space $V_C$ coordinated with q-Lie bracket, namely, that
$a^*=-\overline a$ for any $a\in V_C$ and the following one $u^*=u,\,
v^*=-v,\,w^*=w$. \end{proposition}
{\bf Proof.} Choose a decomposition of $u^{*},v^{*},w^{*}$ over the base
$$ u^{*} = \alpha_{1}u + \beta_{1}v + \gamma_{1}w $$
$$ v^{*} = \alpha_{2}u + \beta_{2}v + \gamma_{2}w $$
$$ w^{*} = \alpha_{3}u + \beta_{3}v + \gamma_{3}w, $$
where $\alpha_{i},\,\beta_{i},\,\gamma_{i},\, i=1,2,3$ are complex
coefficients.
 We want to find them in the accordance with the compatibility condition (1).

It is easy to see that the relation
$$[u^{*},\,u^{*}] = - [u,\,u]^{*} =0$$
implies $\beta_{1}=0$.  Analogicaly, from
$$[w^{*},\,w^{*}] = -  [w,\,w]^{*} = 0$$
we get $\beta_{3} = 0$.

 From the relation
$$[v^{*},\,v^{*}] = - [v,\,v]^{*} = -(1 - q^{2})Mv^{*}$$
we deduce that
$\beta_{2}^{2} + \beta_{2} = 0$, i.e. $\beta_{2}=0$ or $\beta_{2}=-1$.

The relation
$$[w^{*},\,u^{*}] = - [w,\,u]^{*} = (q + q^{-1})^{-1}v^{*}$$
implies $\alpha_2=\gamma_2=0$. If $\beta_{2}=0$ then $v^*=0$, hence
$\beta_{2}=-1$.

At last from
$$[u^{*},\,v^{*}] = - [u,\,v]^{*} = q^{2}Mu^{*}\,\, {\rm and} \,\,
[w^{*},\,v^{*}] = - [w,\,v]^{*} = -Mw^{*}$$
we get $\gamma_{1}=\alpha_{3}=0$.

Thus, we have $v^*=-v,\, u^*=\alpha_1 u,\, w^*=\gamma_3 w$. It is easy to see
that
only two cases are possible $\alpha_1=\gamma_3=-1$ and $\alpha_1=\gamma_3=1$.
This yields the statement.

Finaly we can extend the above involutions to the algebra $\ahqc$ by means
of the relation $*\mu=\mu(*\ot *)\tS$. We leave to the reader to check
the fact that this way to extend the involutions is not contradictory.

Let us note that this proposition implies that involution $u^*=w, v^*=v$ which
corresponds to the compact form of the Lie algebra $sl(2,\,C)$ in the classical
case is not allowed in the braided case and moreover it can not be deformed
in the sense above.

\section{\bf Acknowledgements.} The authors want to thank Y.Kosmann-Schwarzbach
for valuable remarks on the text and to B.Enriquez for helpful discussions.
V.R. is greatly aknowledged the hospitality of Centre de Math\'ematique de
l'Ecole Polytechnique and of the Institute of Theoretical Physics of Uppsala
University.
His work was partially supported by CNRS and by RFFR-MF-95.

 \end{document}